\documentclass{amsart}

\newtheorem{teo}{Theorem}[section]
\newtheorem{lema}{Lemma}[section]

\newtheorem{defi}{Definition}[section]
\newtheorem{cor}{Corollary}[section]
\newtheorem{rem}{Remark}[section]

\def\dem{\noindent{\bf Proof: }}

\def\ve{\varepsilon}
\def\RR{{\mathbb{R}}}
\def\NN{{\mathbb{N}}}

\usepackage{amsfonts}

\begin{document}

\title[Adaptive numerical method]
{Adapting the time step to recover the asymptotic behavior in a
blow-up problem}
\author{Pablo Groisman}
\thanks{Partially supported by Universidad de Buenos Aires under
grant TX48, by ANPCyT PICT No. 03-00000-05009. J. D. Rossi is a
member of CONICET.} \keywords{Blow-up, parabolic equations,
asymptotic behavior.\\
\indent 2000 {\it Mathematics Subjet Classification.} 65M60,
65M20, 35K60, 35B40.}
\address{Departamento de Matem\'atica, FCEyN, UBA \\ (1428)
Buenos Aires, Argentina. } \email{pgroisma@dm.uba.ar}

\begin{abstract}
The equation $u_t = \Delta u + u^p$  with homegeneous Dirichlet
boundary conditions has solutions with blow-up if $p>1$. An
adaptive time-step procedure is given to reproduce the asymptotic
behvior of the solutions in the numerical approximations. We prove
that the numerical method reproduces the blow-up cases, the
blow-up rate and the blow-up time. We also localize the numerical
blow-up set.
\end{abstract}
\maketitle

\section{Introduction.}

\setcounter{equation}{0}

We study the behavior of an adaptive time step
procedure for the following parabolic problem
\begin{equation}
\left\{
\begin{array}{lr}
u_t   =  \Delta u  + u^p \qquad  & \mbox{ in } \Omega \times
[0,T), \\ u (x,t)  =  0\qquad  & \mbox{ on }\partial \Omega \times
[0,T), \\  u(x,0)  =  u_0(x) > 0\qquad  & \mbox{ on } \Omega.
\end{array}
\right.   \label{1.1}
\end{equation}
Where $p$ is superlinear ($p>1$) in order to have solutions with
blow-up.  We assume $u_0$ is regular and $\Omega \subset \RR^d$ is
a bounded smooth domain in order to guarantee that $u \in
C^{2,1}$. A remarkable fact in this problem is that the solution
may develop singularities in finite time, no matter how smooth
$u_0$ is. For many differential equations or systems the solutions
can become unbounded in finite time (a phenomena that is known as
blow-up). Typical examples where this happens are problems
involving reaction terms in the equation like (\ref{1.1}) where a
reaction term of power type is present and so this blow up
phenomenum occurs in the sense that there exists a finite time $T$
such that $\lim_{t\to T} \| u(\cdot,t)\|_\infty=+\infty $ for
initial data large enough (see \cite{P} \cite{SGKM} and the
references therein). The blow-up set, which is defined as the set
composed of points $x \in \Omega$ such that $u(x,t) \to +\infty$
as $t \to T$, is localized in small portions of $\Omega$, in
\cite{V} is proved that the $(d-1)$ dimensional Hausdorff measure
of the blow-up set is finite. The blow-up rate at these points is
given by $u (x, t) \sim (T-t)^{-\frac{1}{p-1}}$, moreover

$$\lim_{t\to T}
(T-t)^{\frac{1}{p-1}}\|u(\cdot,t)\|_{L^{\infty}(\Omega)}=C_{p},
\qquad C_{p}=\left ( \frac{1}{p-1} \right )^{\frac{1}{p-1}}$$
 (see \cite{GK1},\cite{GK2}).

We remark that these results hold if $p$ is subcritical
($p<\frac{d+2}{d-2}$ if $d \ge 3$). For supercritical $p$ the
solutions may present different behaviors. For that reason we
assume $p$ is subcritical along the paper, however the asymptotic
properties of the numerical schemes described above hold for every
$p>1$. This is a difference between the continuous solutions and
their approximations.

Since the solution $u$ develops a singularity in finite time, we
investigate what can be said about the asymptotic behavior (close
to the blow-up time) of numerical approximations for this kind of
problems. In \cite{GR} the authors analyze a semidiscrete scheme
(keeping $t$ continuous) in an interval, $\Omega=(0,1)$. Here we
generalize that analysis for several space dimensions an introduce
the adaptive discretization in time.


For previous work on numerical approximations of (\ref{1.1}) we
refer to \cite{ALM1}, \cite{ALM2}, \cite{BK}, \cite{BHR},
\cite{C}, \cite{LR}, the survey \cite{BB} and references therein.

As a first step to introduce the method we propose a method of
lines, that is: we discretize the space variable, keeping $t$
continuous. In this stage we consider a general method with
adequate assumptions. More precisely, we assume that for every
$h>0$ small ($h$ is the parameter of the method), there exists a
set of nodes $\{x_1,\dots, x_N\}\subset \overline{\Omega}$
($N=N(h)$), such that the numerical approximation of $u$ at the
nodes $x_k$, is given by
$$U(t)=( u_1(t),\dots,u_{N}(t)).$$
That is $u_k(t)$ stands for an approximation of $u(x_k,t)$. We
assume that $U$ is the solution of the following ODE
\begin{equation}
\label{semidiscrete} M U'(t) = - A U (t) + MU^{p},
\end{equation}
with initial data given by $u_k(0) =  u_0(x_k)$. In
(\ref{semidiscrete}) and hereafter, all operations between vectors
are understood coordinate by coordinate.

The precise assumptions on the matrices involved in the method
are:
\begin{enumerate}
\item[(P1)] $M$ is a diagonal matrix with positive entries $m_k$.
\item[(P2)] $A$ is a nonnegative symmetric matrix, with
nonpositive coefficients off the diagonal (i.e. $a_{ij} \le 0 $ if
$i\ne j$) and $a_{ii}>0$. \item[(P3)] $\sum_{k=1}^N a_{ik} \ge 0$.

\end{enumerate}

As an example, we can consider a linear finite element
approximation of problem~(\ref{1.1}) on a regular acute
triangulation of $\Omega$ (see \cite{Ci}). In this case, let $V_h$
be the subspace of $H_0^1(\Omega)$ consisting of piecewise linear
functions on the triangulation. We impose that the finite element
approximation $u_h: [0,T_h) \to V_h$ verifies for each
$t\in[0,T_h)$
$$
\int_{\Omega} ((u_h)_t v)^I = - \int_\Omega \nabla u_h \nabla v +
\int_\Omega ((u_h)^p v)^I,
$$
for every $v \in V_h$. Here $(\cdot)^I$ stands for the linear
Lagrange interpolate at the nodes of the mesh. These conditions
imply that the vector $U(t)$, the values of $u_h(\cdot,t)$ at the
nodes $x_k$, must verify a system of the form
(\ref{semidiscrete}). In this case $M$ is the lumped mass matrix
and $A$ is the stiffness matrix. The assumptions on the matrices
$M$ and $A$ hold as we are considering an acute regular mesh. We
observe that in this case $u_h=U^I$.

As another example, if $\Omega$ is a cube, $\Omega = (0,1)^d$, we
can use a semidiscrete finite differences method to approximate
the solution $u(x,t)$ obtaining an ODE system of the form
(\ref{semidiscrete}).

We also need some kind of convergence result for the scheme, we
will state the precise hypotheses concerning convergence in the
statement of each theorem. Finally, in the Appendix we prove an
$L^\infty $ convergence theorem under the consistency assumption
an give some examples. Now we state the two convergence hypotheses
that we may need.

\begin{enumerate}
\item [(H1)] For every $\tau >0$ $\| u - u_h
\|_{L^\infty(\Omega\times [0,T-\tau])} \to 0$ as $h \to 0$\\

\item [(H2)] $\| u - u_h \|_{H^1_0(\Omega)}(t) \to 0$ as $h \to 0$
for a.e. $t$
\end{enumerate}

Writing these equations explicitly we obtain the following ODE
system,
\begin{equation}\label{system}
m_{k} u_k'(t) = -\sum_{i=1}^N a_{ki} u_{i}(t) + m_{k}u_k^p(t),\  1\le
k\le N,
\end{equation}
with initial data $u_k (0)  = u_0(x_k)$.

In \cite{GR} the authors studied these kind of schemes in one
space dimension, obtaining results similar to those stated below.

Once obtained the ODE system, the next step is to discretize the
time variable $t$. In \cite{ADR} the authors suggest an adaptive
time step procedure to deal with the heat equation with a
nonlinear flux boundary condition. They analyze explicit Euler and
Runge-Kutta methods, however all these methods have to deal with
restrictions in the time-step. In this work we first analyze an
explicit Euler method and next we introduce an implicit scheme in
order to avoid the time-step restrictions. We use $U^j=(u_1^j,
\dots , u_N^j)$ for the values of the numerical approximation at
time $t_j$, and $\tau_j=t_{j+1}-t_{j}$. When we consider the
explicit scheme, $U^j$ is the solution of
\begin{equation}  \label{ec.explicit}
\begin{array}{l}
M U^{j+1}= M U^j + \tau_j \left (- AU^j + M (U^j)^p\right ) \\
U(0)=u_0^I,
\end{array}
\end{equation}
or equivalently, if we denote $\partial u_i^{j+1}=
\frac{1}{\tau_j}(u_i^{j+1} - u_i^j)$

\begin{equation}
\begin{array}{l}
\label{td} m_i\partial u_i^{j+1} =\displaystyle -\sum_{k=1}^N
a_{ki} u_{i}^j +
m_{i}(u^j_i)^p,\  1\le i\le N \\
u^0_i =  u_0(x_i), \qquad 1\le
i\le N+1.
\end{array}
\end{equation}
While for the implicit scheme $U^j$ is the solution of
\begin{equation}  \label{ec.implicit}
\begin{array}{l}
M U^{j+1}= M U^j - \tau_j \left ( AU^{j+1} + M (U^j)^p\right ) \\
U(0)=u_0^I.
\end{array}
\end{equation}

Note that the scheme is not totally implicit since the nonlinear
source $u^p$ is evaluated at time $t^j$ while the discrete
laplacian ($A$) is evaluated a time $t^{j+1}$. This mixture makes
the scheme free of time-step restrictions while the explicit
evaluation of $(U^j)^p$ avoids the problem of solving a nonlinear
system in each step.

Now we choose the time steps $\tau_j=t_{j+1}-t_j$ in such a way
that the asymptotic behavior of the discrete problem is similar to
the continuous one. We will fix $\lambda$ small and take

$$\tau_j=\frac{\lambda}{ (w^j)^p},$$
where
$$w^j= \sum_{k=1}^N m_ku^j_k$$

This choice for the time step is inspired by \cite{ADR}, where the
authors use an adaptive procedure similar to this. They adapt the
time step in a similar way but using the maximum ($L^\infty$-norm)
instead of $w^j$ ($L^1$-norm). In their problem the maximum is
fixed at the right boundary node (i.e. $\|U^j\|_\infty =
u^j_{N+1}$). In this problem, the maximum (the node $k$ such that
$u_k^j= \|U^j\|_\infty$) can move from one node to another as $t$
goes forward. So it is better to deal with $w^j$ since, for
example, we can compute its derivative. Anyway, as we use a fixed
mesh for the space discretization we can compare both norms. The
motivation of this choice for the time-step is that, as will be
shown, the behavior of $w^j$ is given by
$$
\partial w^{j+1} \sim (w^j)^p.
$$
Hence, with our selection of $\tau_j$ we can obtain
$$
w^{j+1} \sim w^j + \tau_j (w^j)^p = w^j + \lambda \sim w^0 +
(j+1)\lambda,
$$
and we obtain the asymptotic behavior of $w^j$, which is, as we
will see, similar to the one for the continuous solution.

When we deal with the explicit scheme we will also require

\begin{equation}
\label{restriction}
 \lambda < \min_{1\le i \le
N}\frac{m_i}{a_{ii}}(w^0)^p.
\end{equation}
Then we study the asymptotic properties of the numerical schemes.
We will say that a solution of \eqref{ec.explicit} (or
\eqref{ec.implicit}) blows up if

$$\lim_{j\to \infty}\|U^j\|_\infty = \infty, \qquad \mbox{ and }
\qquad T_{h,\lambda}:=\sum_{j=1}^\infty \tau_j < \infty,$$ we call
$T_{h,\lambda}$ the blow-up time. To describe when the blow-up
phenomena occurs in the discrete problem we introduce the
following functional $\Phi_h:~\RR^N~\to~\RR$.

$$\Phi_{h} (U) \equiv \langle A U,  U \rangle
 - \langle \frac{1}{p+1}MU^{p+1}, ME  \rangle, $$
 where $E=(1,1,\dots,1)$. The functional $\Phi_h$ is a discrete version of

$$\Phi (u) (t) \equiv \int_\Omega \frac{|\nabla u (s,t) |^2  }{2} \ ds-
 \int_\Omega \frac{ (u(s,t))^{p+1} }{p+1} \ ds .$$

This functional characterize the solutions with blow-up in the
continuous problem: in \cite{CPE}, \cite{GK2} it is proved that
$u$ blows up at time $T$ if and only if $\Phi(u)(t) \to -\infty$
as $t \to T$. We prove a similar result for the discrete
functional $\Phi_h$ and this allows as to prove that if the
continuous solution has finite time blow-up its numerical
approximation also blows up when the parameters of the method are
small enough.

Next we study the asymptotic behavior for the numerical
approximations of the solutions with blow-up and we find that they
behave very similar to the continuous ones. In fact we find that
if $u_{h,\lambda}$ is a numerical solution with blow-up at time
$T_{h,\lambda}$ its behavior is given by

$$\displaystyle \max_{1 \le i \le N} u_i^j \sim (T_{h,\lambda} -
t^j)^{-1/(p-1)}.$$

We use the notation $f(j) \sim g(j)$ to mean that there exist
constants $c, C
>0$ independent of $j$ (but they may depend on $h$) such that

$$ cg(j) \le f(j) \le Cg(j)$$

Moreover, the numerical schemes recover the constant $C_p$ in the
sense that

$$ \displaystyle \lim _{j \to \infty} \max_{1 \le i \le N} u_i^j
(T_{h,\lambda} - t^j)^{1/(p-1)} = C_p .$$

The functional $\Phi_h$ is also useful to prove convergence of the
blow-up times. Unfortunately we can only prove the convergence of
an iterated limit,
$$\lim_{h\to 0}\lim_{\lambda \to 0} T_{h,\lambda} = T.$$

By means of the numerical blow-up rate we observe a propagation
property for blow-up points. We prove that the nodes adjacent to
those that blow-up as the maximum may also blow-up (opposite to
the continuous problem), but they did it with a slower rate and
the number of adjacent blow-up nodes is determined only by $p$ and
is independent of $h$ and $\lambda$.

In other words, if we call $B^*(U)$ the set of nodes $k$ such that
$$ \displaystyle \lim _{j \to \infty} u_k^j
(T_{h,\lambda} - t^j)^{1/(p-1)} = C_p ,$$ the number of blow-up
points outside $B^*(U)$ depends (explicitly) only on the power
$p$.

We split the paper in two parts, in the first part we deal with
the explicit scheme and in the second one we develop the analysis
for the implicit method.

\section{The explicit scheme}
\setcounter{equation}{0}
 The main tool in our proofs is a
comparison argument, so first of all we prove a lemma which states
that this comparison argument holds. Since we need restrictions in
the time-step to prove this lemma they are essential for every
result in this section. That is not the case of the implicit
scheme.

\begin{defi}
We say that $(Z^j)$ is a supersolution (resp.: subsolution) for
$\eqref{ec.explicit}$ if verifies the equation with an inequality
$\ge$ ($\le $) instead of an equality.
\end{defi}

\begin{lema}
\label{pmaximo} Assume the time step verifies
$$\tau_j < \min_{1\le i \le N} \frac{m_i}{a_{ii}}.$$
Let $(\overline{U}^j), (\underline{U}^j)$ a super and a
subsolution respectively for $\eqref{ec.explicit}$ such that
$\underline{U}^0<\overline{U}^0$, then
$\underline{U}^j<\overline{U}^j$ for every $j$.
\end{lema}

\dem Let $Z^j=\overline{U}^j - \underline{U}^j$, by an
approximation argument we can assume that we have strict
inequalities in \eqref{ec.explicit}, then  $(Z^j)$ verifies

\begin{equation}
\label{supersol}
\begin{array}{lcl}
 M \partial Z^{j+1} & > & -A Z^j +
M((\overline{U}^j)^p - (\underline{U}^j)^p, \\
Z^0 & > & 0. \end{array}
\end{equation}

If the statement of the Lemma is false, then there exists a first
time $t^{j+1}$ and a node $x_i$ such that $z_i^{j+1} \le 0$. At
that time we have

$$
z_i^{j+1} > (1 - \tau_j \frac{a_{ii}}{m_i} )z_i^j  + \tau_j \left
( -\sum_{k\ne i} a_{ik}z_k^j + (\overline{u}_i^j)^p -
(\underline{u}_i^j)^p \right )\ge 0,$$ a contradiction. \qed

\subsection{When does the solution blow up.} In this section we
find conditions under which the solution of \eqref{td} blows up,
we begin with some lemmas.

As the matrix $A$ is a symmetric (property (P2)), there exists a
basis of eigenvectors for the following eigenvalue problem
$$A\phi_i=\lambda_iM\phi_i.$$
We call $\eta=\eta(h)$ the greatest eigenvalue of this problem,
that is
$$0 \le \lambda_i \le \eta(h).$$

\begin{lema}
\label{comp} For every $y \in \RR^N$ there holds
$$ \langle Ay,y \rangle \le \eta(h) \langle My,y \rangle.$$
\end{lema}

\dem As the matrix $M$ defines a scalar product in $\RR^N$, we can
assume that the eigenvectors $\phi_i$ are normalized such that
$$\langle M\phi_i, \phi_j \rangle = \delta_{ij}.$$
Let $y \in \RR^N, y= \sum_{i=1}^N \alpha_i \phi_i$, then
\begin{eqnarray}
\nonumber \langle Ay, y \rangle & = & \left \langle \sum_{i=1}^N
\alpha_i \lambda_i M \phi_i, \sum_{j=1}^N \alpha_j \phi_j \right
\rangle\\
\nonumber & = & \sum_{i=1}^N \alpha^2_i \lambda_i \langle M\phi_i,
\phi_i \rangle\\
\nonumber & \le & \eta(h)\langle My, y \rangle.
\end{eqnarray}
\qed

\begin{lema}
\label{siesgrandeexplota} If  $(U^j)_{j\ge 0}$ is large enough,
then blows up in finite time.
\end{lema}

\dem Recall the definition of

$$w^j=  \sum_{i=1}^N m_i u^j_i,$$
and observe that there exists constants $c,C>0$ that depend only
on $h$ such that
$$ c\max_{1\le k \le N} u_i^j \le w^j \le C\max_{1\le i \le N}
u_i^j.$$

Now we observe that we can choose $j \ge j_0$ in order to get
$w^j$ large enough to verify
\begin{eqnarray}
\nonumber w^{j+1} & = & w^j - \tau_j \sum_{k=1}^N\sum_{i=1}^N
a_{ik}u_k^j + \tau_j
\sum_{i=1}^{N}m_i (u_i^j)^p\\
\nonumber & \ge & w^j - \frac12 \tau_j \sum_{i=1}^N m_{i}(u_i^j)^p
+ \tau_j
\sum_{i=1}^{N}m_i (u_i^j)^p\\
\nonumber & \ge & w^j + \frac{\tau_j}{2} \sum_{i=1}^N
m_{i}(u_i^j)^p \\
\nonumber & \ge  &  w^j   + c \tau_j (w^j)^p\\
\nonumber & = & w^j + c\lambda.\\
\nonumber
\end{eqnarray}

Applying this inequality inductively we obtain
$$ w^j \ge w^{j_0} + c\lambda(j-j_0) \ge cj,$$
hence  $w^j \to \infty$ as $j \to \infty$, it remains to check
that $\sum \tau_j < \infty$, to do that we observe
$$\tau_j = \frac{\lambda}{(w^j)^p} \le \frac{\lambda}{(w^0 +
jc\lambda)^p},$$ and
$$ \sum_{j=1}^\infty \frac{\lambda}{(w^0 +
j(c\lambda))^p} \le \int_0^\infty \frac{\lambda}{(w^0 +
cs\lambda)^p} \, ds < \infty.$$ This completes the proof. \qed

{\rem In the course of the proof of the above Lemma we showed not
just that $w^j$ blows up, we also proved $w^j \ge cj$.}

Now we are going to prove the reverse inequality to obtain the
asymptotic behavior of $\|U^j\|_\infty$.

\begin{lema} If $(U^j)$ is unbounded then
$$\|U^j\|_\infty \sim w^j \sim j$$
\end{lema}

\dem The relation $\|U^j\|_\infty \sim w^j$ is trivial since they
define equivalent norms in $\RR^N$. So we just have to prove $w^j
\le Cj$. To do that we observe that

\begin{eqnarray}
\nonumber w^{j+1} & = & w^j - \tau_j \sum_{k=1}^N\sum_{i=1}^N
a_{ik}u_k^j + \tau_j
\sum_{i=1}^{N}m_i (u_i^j)^p\\
\nonumber & \le & w^j  + \tau_j \sum_{i=1}^{N} m_i (u_i^j)^p\\
\nonumber & \le  &  w^j   + C \tau_j (w^j)^p\\
\nonumber & = & w^j + C\lambda.\\
\nonumber
\end{eqnarray}

We apply this inequality inductively again to get
$$ w^j \le w^{0} + C\lambda j \le Cj,$$
as we wanted to prove. \qed

\medskip

\begin{teo}
\label{funcional} Assume the time step $\tau_j$ verifies the
restriction
\begin{equation}
\label{restriccion.tiempo} \tau_j < \frac{2}{p(w^{j+1})^{p-1} +
\eta(h)}.
\end{equation}
Then positive solutions of \eqref{td} blow up if there exists
$j_0$ such that $\Phi_{h}(U^{j_0})<0$.
\end{teo}
We remark that the condition $\Phi_{h}(U^{j_0})<0$ is similar to
the one for the blow-up phenomena in the continuous problem, in
fact for the continuous problem this is a necessary condition.


\dem First we observe that $\Phi_{h} (U^j)$ decreases with $j$, in
order to do that we take inner product of \eqref{ec.explicit} with
$U^{j+1}-U^j$  to obtain
\begin{eqnarray*}
\nonumber  0 & = & \langle \frac{1}{\tau_j} M(U^{j+1} -  U^j)  +
AU^j
- M(U^j)^p,  U^{j+1} -  U^j\rangle\\
\nonumber &=&\tau_j \langle M\partial U^{j+1}, \partial U^{j+1}
\rangle + \Phi_h(U^{j+1}) - \Phi_h(U^j)- \frac12 \langle AU^{j+1},
U^{j+1} \rangle\\
& + &\langle AU^j, U^{j+1}\rangle - \frac12 \langle AU^j,
U^j\rangle - \frac12 \langle Mp(\xi^j)^{p-1}, (U^{j+1} - U^j)^2
\rangle.
\end{eqnarray*}
Hence we obtain,
$$
\begin{array}{l}
\Phi_h(U^{j+1})  - \Phi_h(U^j)\\
\\
\quad \displaystyle \le \tau_j(\tau_j\frac{p(w^{j+1})^{p-1}}{2}
-1) \langle M\partial U^{j+1},
\partial U^{j+1} \rangle +  \, \frac{\tau_j^2}{2}\langle A\partial
U^{j+1}, \partial U^{j+1} \rangle\\
\\
\quad \displaystyle\le  \tau_j(\tau_j \displaystyle
\frac{p(w^{j+1})^{p-1}}{2} + \frac{\eta(h)\tau_j}{2} -1)\langle
M\partial U^{j+1},
\partial U^{j+1} \rangle  \le  0,
\end{array}
$$
due the restriction in the time step $\tau_j$ and Lemma
\ref{comp}. Actually $\Phi_h(U^{j+1}) < \Phi_h(U^j)$ unless
$(U^j)$ is independent of $j$. So, $\Phi_h$ is a Lyapunov
functional for \eqref{ec.explicit}. Next we observe that the
steady states of \eqref{ec.explicit} have positive energy. Let
$(W^j)=W$ be a stationary solution of \eqref{ec.explicit}, then
$$0=-AW + MW^p.$$

Multiplying by $W/2$ we obtain,
\begin{eqnarray}
\nonumber 0&=& -\frac12 \langle AW, W \rangle + \frac{p+1}{2}
\frac{1}{p+1} \langle MW^{p},W\rangle\\
\nonumber & \ge & -\Phi_h(W).
\end{eqnarray}

Now, assume $(U^j)$ is a bounded solution of \eqref{ec.explicit},
then there exists a convergent subsequence of $(U^j)$ that we
still denote $(U^j)$. Its limit $W$ is a steady state with
positive energy.

As $\Phi_{h,\lambda} (U^j)$ decreases and there exists $j_0$ such
that $\Phi_{h,\lambda} (U^{j_0}) < 0$ then $\Phi_{h,\lambda} (W) <
0$, a contradiction. We conclude that $(U^j)$ is unbounded and by
Lemma \ref{siesgrandeexplota} has finite time blow-up. \qed

\begin{cor}
\label{cor1.e} Assume the time-step restriction of the above
theorem and the convergence hypotheses {\rm (H1), (H2)}. Let $u_0$
an initial data for \eqref{1.1} such that $u$ blows up in finite
time $T$. Then $u_{h, \lambda}$ blows up in finite time $T_{h,
\lambda}$ for every $h$, $\lambda = \lambda(h)$ small enough.
Moreover

$$\lim_{h \to 0} \lim_{\lambda \to 0} T_{h, \lambda} = T$$

\end{cor}

{\rem If the fully-discrete method converges in $H_0^1(\Omega)$
a.e. $t$ then $\lambda$ can be chosen independent of $h$.}

\dem We observe that if $u$ blows up in finite time $T$ then (see
\cite{CPE}, \cite{GK2})

$$ \Phi (u) (t) \equiv \int_\Omega \frac{|\nabla u (s,t) |^2  }{2} \
ds-
 \int_\Omega \frac{ (u(s,t))^{p+1} }{p+1} \ ds \to -\infty \qquad (t
 \nearrow T).$$

This implies that there exists a time $t_0<T$ with
$\Phi(u)(t_0)<0$. Now we use the convergence of $u_{h,\lambda}$ to
$u_h$ in $[0,t_0]$ and the convergence hypothesis (H1) to see that
$$
\lim_{h\to 0}\lim_{\lambda \to 0}
\Phi_{h,\lambda}(u_{h,\lambda})(t_0) = \Phi(u)(t_0).
$$
So for $h,\lambda=\lambda(h)$ small enough we get
$\Phi_{h,\lambda}(u_{h,\lambda})(t_0)<0$ and so, by the above
Theorem $(U^j)$ blows up.

Now we turn our attention to the blow-up times. In \cite{GR} it is
proved that the blow-up time of the semi-discrete solutions
(solutions of \eqref{semidiscrete}), that we are going to denote
$T_h$, converges, as $h \to 0$, to $T$. In that work the poof is
given for $\Omega=(0,1)$ and a finite element method. Despite for
this case the proof is very similar, we sketch it for the sake of
completeness. In \cite{GR} is proved that if the continuous
solution blows up then for every $h$ small enough the semidiscrete
scheme also does. Hence we can assume that the semidiscrete
solution $U(t)$ is large enough in order to verify

$$\frac{d}{dt} \langle MU(t),U(t) \rangle =
2  \langle MU'(t),U(t) \rangle = $$
$$2 \langle -AU(t),U(t) \rangle
+ 2 \langle MU^p(t),U(t) \rangle =$$
$$ -4\Phi_h(U(t)) +  \frac{2(p-1)}{p+1}
\langle MU^p(t),U(t) \rangle \ge $$
$$ 4 |\Phi_h(U(t))|
+ \frac{2(p-1)}{p+1}(\langle MU(t),U(t) \rangle)^\frac{p+1}{2}.$$
Integrating between $t_0$ and $T_h$ we obtain

\begin{equation}
\label{cota} (T_h -t_0 ) \le \frac{C}{ ( - \Phi_h (U(t_0))
)^{\frac{p-1}{p+1}} }. \end{equation} where $C$ depends only on
$p$.

Given $\ve >0$, we can choose $M$ large enough to ensure that
$$ \left( \frac{ C}{ M^{\frac{p-1}{p+1}} } \right)  \le \frac{\ve}{2}.$$
As $u$ blows up at time $T$ we can choose $\tau<\frac{\ve}{2}$
such that
$$- \Phi(u(\cdot, T-\tau) \ge 2M.$$
If $h$ is small enough,
$$ - \Phi_h ( U(T-\tau))  \ge M,$$
and hence
$$T_h - (T-\tau) \le
\left( \frac{C}{  ( - \Phi_h (U(T - \tau)) )^{\frac{p-1}{p+1}}
}\right) \le \left( \frac{C}{ M^{\frac{p-1}{p+1}} } \right )  \le
\frac{\ve}{2} .$$ Therefore,
$$|T_h-T|\le |T_h-(T-\tau)|+|\tau|<\ve.$$

We have proved $\lim_{h \to 0} T_h = T$, so we just have to prove
that for fixed $h$

$$\lim_{\lambda \to 0} T_{h, \lambda} = T_h.$$
To do that we observe that from Lemma \ref{siesgrandeexplota} we
know that there exists $j_0$, that does not depend on $\lambda$
such that for $j \ge j_0$

$$ w^j \ge w^{j_0} + c\lambda(j-j_0),$$
hence
\begin{eqnarray}
\nonumber T_{h, \lambda} - t^j & = & \sum_{l=j+1}^\infty \tau_l =
\sum_{l=j+1}^\infty \frac{\lambda}{(w^l)^p}\\
\nonumber & \le & \sum_{k=j+1}^\infty \frac{\lambda}{(w^{j_0} +
c\lambda(l-j_0))^p} \le \int_j^\infty \frac{\lambda}{(w^{j_0} +
c\lambda(s-j_0))^p} \, ds\\
\nonumber & = & \frac{1}{c}\int_{w^{j_0} + c\lambda(j-j_0)}^\infty
\frac{1}{s^p} \, ds \le \frac{1}{c}\int_{w^{j_0}}^\infty
\frac{1}{s^p} \, ds.\\
\end{eqnarray}
This holds for any $j_0$ large enough and for every $j \ge  j_0$.
In particular we get

$$T_{h, \lambda} - t^j \le \frac{1}{c}\int_{w^{j}}^\infty
\frac{1}{s^p} \, ds.$$

This inequality has a great meaning. It says that if $w^j$ is
large, then $t^j$ is close to blow up (independent of $\lambda$).
So now, given $\ve >0$ we can choose $K$ large enough in order to
have

$$\frac{1}{c}\int_{K}^\infty \frac{1}{s^p} \, ds \le
\frac{\ve}{3}, \qquad K^{-p} < \frac{\ve}{3}.$$

Next we choose $\tau < \frac{\ve}{3}$ such that $\sum m_k u_k(T_h
-2\tau) \ge 2K$ (remember that $(u_1(t), \dots , u_N(t))$ is the
solution of the semidiscrete scheme). For
$\lambda=\lambda(h,\tau)$ small enough we get, from (H2), that
$w^{j} \ge K$ for every $j$ such that $T_h - 2\tau \le t^j \le T_h
- \tau$. We choose one of those $j$ and compute

\begin{eqnarray}
\nonumber |T_{h, \lambda} - T_h| & \le &|T_{h, \lambda} - t^j| +
|t^j - T_h|\\
\nonumber & \le & \frac{1}{c}\int_{K}^\infty \frac{1}{s^p} \, ds +
2\tau\\
\nonumber & \le & \ve
\end{eqnarray}
\qed

\subsection{Blow-up rate}

In this section we study the asymptotic behavior of numerical
solutions with blow-up.

\begin{teo}
\label{tasa.e}
 Let $u_{h, \lambda}$ a solution with blow-up at time
$T_{h,\lambda}$, then

$$\displaystyle \max_{1 \le i \le N} u_i^j \sim (T_{h,\lambda} -
t^j)^{-1/(p-1)}.$$ Moreover
$$ \displaystyle \lim _{j \to \infty} \max_{1 \le i \le N} u_i^j
(T_{h,\lambda} - t^j)^{1/(p-1)} = C_p = \left (
\frac{1}{p-1}\right )^{1/(p-1)}.$$

\end{teo}

We want to remark that this is the behavior of the continuous
solutions with blow-up.

\dem We know from Lemma \ref{siesgrandeexplota} that $w^j = \sum
m_iu_i^j$ verifies
$$w^{j+1} \ge w^j + c\tau_j(w^j)^p,$$
so we have

\begin{eqnarray}
\nonumber (T_{h,\lambda} - t^j) & = & \sum_{k=j+1}^\infty \tau_j =
\sum_{k=j+1}^\infty \frac{\lambda}{(w^j)^p} \\
\nonumber & \le & \int_j^\infty \frac{\lambda}{(w(s))^p}\, ds.\\
\end{eqnarray}
Here $w(s)$ is the linear interpolant of $(w^j=w(j))$, hence for
$j \le s \le j+1$ we have $w'(s)= w^{j+1}-w^j \ge c\lambda$, and
so
$$
\int_j^\infty \frac{\lambda}{(w(s))^p}\, ds  \le \int_{w^j}^\infty
\frac{\lambda}{cv^p\lambda} \, dv \le \frac{1}{c(p-1)}\left (
\frac{1}{w^j} \right )^{p-1},$$ or equivalently
$$ \max_{1\le i \le N} u_i^j \le Cw^j \le C(T_{h,\lambda} -
t^j)^{-1/(p-1)}.$$ The inverse inequalities can be handled in a
similar way to obtain
$$ \max_{1\le i \le N} u_i^j \sim w^j \sim (T_{h,\lambda} -
t^j)^{-1/(p-1)}.$$

Next we recover the constant $C_p$ in the asymptotic behavior of
the numerical solution, to do that we will change variables but
first we prove a short lemma.

\begin{lema}
\label{limite}
$$\lim_{j \to \infty}\frac{T_{h, \lambda} - t^j}{T_{h,
\lambda}-t^{j+1}}=1$$
\end{lema}

\begin{dem}
$$ \displaystyle
1 \le \frac{T_{h, \lambda} - t^j}{T_{h, \lambda}-t^{j+1}} =
\frac{\sum_{k=j+1}^\infty \tau_k}{\sum_{k=j+2}^\infty \tau_k} = 1
+ \frac{ \tau_{j+1}}{\sum_{k=j+2}^\infty \tau_k}$$ $$ \le 1+
\frac{\lambda/(w^{j+1})^p}{C/ (w^{j+1})^{p-1}} \to 1.
$$\qed
\end{dem}

Now we change variables, in a way inspired by
\cite{GK1},\cite{GR}. Let $(Y^j)$ be defined by

\begin{equation}
y_i^j=u_i^j(T_{h,\lambda} - t^j)^{1/(p-1)} \qquad 1 \le i \le N.
\end{equation}

In the sequel of the proof we will use $\Delta y_i^{j+1}$ to
denote

$$\frac{y_i^{j+1} - y_i^j}{\tau_j/(T_{h, \lambda}-t^j)},$$

This can be thought as $\tau_j/(T_{h,\lambda}- t^j)$ to be the
time step in the new variables. With this notation the new
variables verify

\begin{equation}
\begin{array}{rcl}
\label{autosimilares} m_i\Delta y_i^{j+1} &=& \displaystyle
-\frac{(T_{h,\lambda}-t^{j+1})^{\frac{1}{p-1}}}{(T_{h,\lambda}-t^{j})^{\frac{1}{p-1}}}
(T_{h,\lambda}-t^{j})
 \sum_{i=1}^N a_{ki} y_{i}^j\\
 \\
 & & \displaystyle +  m_i
\frac{(T_{h,\lambda}-t^{j+1})^\frac{1}{p-1}}{(T_{h,\lambda}-t^{j})^\frac{1}{p-1}}(y_i^j)^p\\
\\
  &  &+  \displaystyle \frac{(T_{h,\lambda} - t^j) m_i
u_i^{j}}{\tau_j}
 \left( (T_{h,\lambda}-t^{j+1})^{\frac{1}{p-1}}-
(T_{h,\lambda}-t^{j})^{\frac{1}{p-1}} \right ), \\
 \\
y^0_i& = & (T_{h,\lambda})^{1/(p-1)}u_0(x_i), \qquad 1\le i\le
N+1.
\end{array}
\end{equation}

Now assume there exists $\ve >0$ and a subsequence that we still
denote $(y_i^j)$ such that $y_i^j > C_p + \ve$ for some $i=i(j)$.
Then for those $y_i^j$ we have

$$(y_i^j)^p - \frac{1}{p-1}y_i^j > \frac{3 \delta}{m_i}.$$

We also know from the blow-up rate that the new variables $y_i^j$
are bounded and so, applying Lemma \ref{limite} we obtain for  $j$
large enough

\begin{eqnarray}
\label{Cp.e} \nonumber m_i\Delta y_i^{j+1} &\ge& -\delta  + m_i
\left ((y_i^j)^p - \frac{1}{p-1}y_i^j \right )\\
\\
  & & +  \displaystyle  m_i (y_i^j)^p \left (
\frac{(T_{h,\lambda}-t^{j+1})^\frac{1}{p-1}}{(T_{h,\lambda}-t^{j})^\frac{1}{p-1}}
-1 \right )\\
\nonumber & \ge & \delta.
 \\
\end{eqnarray}

This means that actually $y_i^j > C_p + \ve$ for every $j$ large
enough and consequently \eqref{Cp.e} is verified for all those
$j$. So $y_i^j$ is unbounded, a contradiction.

If we assume $y_i^j < C_p - \ve$ arguing along the same lines we
obtain that $y_i^j$ verifies

\begin{eqnarray}
\nonumber m_i\Delta y_i^{j+1} &\le& \delta  + m_i
\frac{(T_{h,\lambda}-t^{j+1})^{\frac{1}{p-1}}}{(T_{h,\lambda}-t^{j})^{\frac{1}{p-1}}}\left
((y_i^j)^p  \frac{1}{p-1}y_i^j \right )\\
\nonumber &  &+  \displaystyle  \frac{m_i}{p-1} y_i^j \left (
\frac{(T_{h,\lambda}-t^{j+1})^\frac{1}{p-1}}{(T_{h,\lambda}-t^{j})^\frac{1}{p-1}}
-
\frac{(T_{h,\lambda}-\xi^{j})^{\frac{1}{p-1}-1}}{(T_{h,\lambda}-t^{j})^{\frac{1}{p-1}-1}}
\right )\\
\nonumber & \le & 2 \delta + C \left ((y_i^j)^p -
\frac{1}{p-1}y_i^j \right ).
\\
\end{eqnarray}

This shows that either $y_i^j \to 0$ as $j \to \infty$ or
$m_i\Delta y_i^{j+1} < - \delta$, which means that $y_i^j$ is not
bounded from below (this is not possible).

We conclude that if $y_i^j$ does not go to zero, then it goes to
$C_p$, as the blow-up rate implies that for every $j$
$$ \max_{1 \le i \le N} y_i^j \ge c,$$
we have
$$\lim_{j \to \infty}\max_{1 \le i \le N} y_i^j =
 \lim_{j \to \infty}\max_{1 \le i \le N} (T_{h, \lambda} -
 t^j)^{1/(p-1)}u_i^j= C_p,
$$
as we wanted to prove. \qed

\subsection{Blow-up set}

Now we turn our attention to the blow-up set. In order to do that
we consider the set $B^*(U)$ composed of those nodes that blow-up
as the maximum (a precise definition of $B^*(U)$ is given in the
introduction) and we study the behavior of the adjacent nodes,
then we repeat the procedure with these nodes.

\begin{defi}\label{dist}
We define the graph with vertices in the nodes and say that two
different nodes are connected if and only if $a_{ij}\ne 0$. We
consider the usual distance between nodes measured as a graph, see
\cite{Ha}. Finally, we denote by $d(k)$ the distance of the node
$x_k$ to $B^*(U)$ also measured as a graph. \end{defi}

We prove that $u_k$ blows up if and only if $d(k)\le K$ where $K$
depends only on $p$,

\begin{teo}
\label{teo5} Let $B^*(U)$ be the set of nodes, $\{x_k\}$, such
that
$$u_k^j \sim (T_{h,\lambda}-t^j)^{- \frac{1}{p-1}}.$$ Then the blow-up
propagates in the following way, let $p>1$ and $K \in \NN_0$ such
that $\frac{K+2}{K+1}<p \le \frac{K+1}K$ ($K$ is the integer part
of $1/(p-1)$). Then the solution of \eqref{ec.explicit} blows up
exactly at $K$ nodes near $B^*(U)$. More precisely,
$$u_k^j \to +\infty \qquad \iff \qquad d(k) \le K.$$ Moreover, if
$d(k)\le K$, the asymptotic behaviour of $(u^j_k)_{j\ge 1}$ is
given by
$$u_{k}^j \sim (T_{h,\lambda}-t^j)^{-\frac{1}{p-1} + d(k) },$$ if
$p \ne \frac{K+1}{K}$ and if $p = \frac{K+1}{K}$, $d(k) =K$
$$u_{k}^j \sim \ln (T_{h,\lambda} -t^j).$$
\end{teo}

\medskip

We want to remark that more than one node can go to infinity, but
the asymptotic behavior imposes $\frac{u_{k}^j}{u_{i}^j}\to 0$ ($j
\to \infty$) if $d(k)>d(i)$.

{\bf Proof of Theorem 1.4} We want to show that the blow-up
propagates $K$ nodes around $B^*(U)$, we begin with a node $x_k$
such that $d(k)=1$. We claim that the behavior of $u_k^j$ is given
by

$$u_k^j \sim \left \{ \begin{array}{lll}
j^{-p+2} & \mbox { if } & p<2\\
\ln{j} & \mbox{ if } & p=2\\
C & \mbox{ if } & p>2,
\end{array} \right.
$$
to prove that we will show that
$$w_{A}^{j}=A\sum_{s=1}^{j}s\tau_{s-1},$$
which has the behavior described above, can be used as super and
subsolution for an equation verified by $u_{k}^{j}$ choosing $A$
appropriately.

We observe that $u_{k}^{j}$ satisfies

$$ m_k\partial u_k^{j+1} =\displaystyle -\sum_{i=1}^N
a_{ik} u_{i}^j + m_{k}(u^j_k)^p $$ $$\sim C_{1}(\max_{1\le i \le
N} u_{i}^{j}) -C_{2}u^{j}_{k} + C_{3}m_{k}(u^j_k)^p.$$ This means
that there exists constants $c_{i}, C_{i}>0, \, i=1,2,3$ such that
for $j$ large enough
\begin{equation}
\label{sub1}
\partial u_k^{j+1} \le C_{1}j  -C_{2}u^{j}_{k} + C_{3}(u^j_k)^p
\end{equation}
and
\begin{equation}
\label{super1}
\partial u_k^{j+1} \ge c_{1}j - c_{2}u^{j}_{k}+ c_{3}(u^j_k)^p
\end{equation}
Now we observe that if $A$ and $j$ are large enough, $w_{A}^{j}$
verifies

\begin{eqnarray}
\nonumber \partial w_{A}^{j} &=& A(j+1)\\
\nonumber & \ge & C_{1}j  - C_{2}w^{j}_{A} + C_{3}(w_{A}^{j})^{p},
\end{eqnarray}
since $(w_{A}^{j})^{p}/j \to 0$ as $j$ goes to infintiy.
Hence $w_{A}^{j}$ is a supersolution for \eqref{sub1} and so
$$u_{k}^{j} \le w^{j}_{A}.$$
On the other hand if we choose $A$ small we get
\begin{eqnarray}\
\nonumber \partial w_{A}^{j} &=& A(j+1) \\
\nonumber & \le & c_{1}j - c_{2}w_{A}^{j} + c_{3}(w_{A}^{j})^{p},
\end{eqnarray}
Hence now we can use $w_{A}^{j}$ as a subsolution for
\eqref{super1} to handle the other inequality. Therefore
$$u_{k}^{j} \sim w_{A}^{j}.$$

We observe that if $p<2$ the node $x_{k}$ is a blow-up node and we
also have the blow-up rate for this node ($u_{k}^{j} \sim
j^{-p+2}$). If $p>2$ this node is bounded. Next we assume $p<2$
(if $p>2$ it is easy to prove that every node $k$ with $d(k) \ge
1$ is bounded) and we are going to find the behavior of a node,
that we still denote $k$, such that $d(k)=2$.That is, it is not
adyacent to $B^*(U)$ and it is adyacent to a node wich has the
behavior $j^{-p+2}$.

Now let
$$w_{A}^{j}=A\sum_{s=1}^{j}\tau_{s}s^{-p+2},$$
and observe that $u_{k}^{j}$ verifies

$$ m_k\partial u_k^{j+1} =\displaystyle -\sum_{i=1}^N
a_{ik} u_{i}^j +
m_{k}(u^j_k)^p \sim C_{1}(j^{-p+2}) -C_{2}u^{j}_{k} +
C_{3}m_{k}(u^j_k)^p.$$
This means that there exists constants $c_{i}, C_{i}>0, \, i=1,2,3$
such that for $j$ large enough
\begin{equation}
\label{sub2}
\partial u_k^{j+1} \le C_{1}j^{-p+2}  -C_{2}u^{j}_{k} + C_{3}(u^j_k)^p
\end{equation}
and
\begin{equation}
\label{super1}
\partial u_k^{j+1} \ge c_{1}j^{-p+2} - c_{2}u^{j}_{k}+ c_{3}(u^j_k)^p
\end{equation}
Now we observe that if $A$ and $j$ are large enough, $w_{A}^{j}$
verifies

\begin{eqnarray}
\nonumber \partial w_{A}^{j} &=& A(j+1)^{-p+2}\\
\nonumber & \ge & C_{1}j^{-p+2}  - C_{2}w^{j}_{A} + C_{3}(w_{A}^{j})^{p},
\end{eqnarray}
since $(w_{A}^{j})^{p}/j^{-p+2} \to 0$ as $j$ goes to infintiy.
Hence $w_{A}^{j}$ is a supersolution for \eqref{sub2} and so
$$u_{k}^{j} \le w^{j}_{A}.$$
On the other hand if we choose $A$ small we get
\begin{eqnarray}\
\nonumber \partial w_{A}^{j} &=& A(j+1)^{-p+2} \\
\nonumber & \le & c_{1}j^{-p+2} - c_{2}w_{A}^{j} + c_{3}(w_{A}^{j})^{p},
\end{eqnarray}
Now we can use $w_{A}^{j}$ as a subsolution for \eqref{super1} to
handle the other inequality. So
$$u_{k}^{j} \sim w_{A}^{j}.$$

We observe that if $p<3/2$ the node $x_{k}$ is a blow-up node and
we also have the blow-up rate for this node ($u_{k}^{j} \sim
j^{-2p+3}$). If $p>3/2$ this node is bounded. In the case $p<3/2$
we repeat this procedure inductively to obtain the theorem. \qed

\section{The implicit scheme}
\setcounter{equation}{0}

In order to avoid the time step restrictions we now introduce an
implicit scheme and prove that similar properties can be observed.
However we have to remark that near the blow-up time the adaptive
procedure forces the time step to be much smaller than the space
discretization parameter $h$. This suggest that an adequate method
could be to begin at time zero with the implicit scheme in order
to avoid time-step restrictions and, as when the solution
increases the time step is reduced, if the solution has finite
time blow-up, then there will be a moment such that the time step
restriction for the explicit scheme will be verified. From that
time one can continue either with the explicit or with the
implicit scheme. As in the explicit scheme section we begin with
the comparison Lemma.

\begin{lema}
\label{pmaximo.i} Let $(\overline{U}^j), (\underline{U}^j)$ a
super and a subsolution respectively for $\eqref{ec.implicit}$
such that $\underline{U}^0<\overline{U}^0$, then
$\underline{U}^j<\overline{U}^j$ for every $j$.
\end{lema}

\dem Let $Z^j=\overline{U}^j - \underline{U}^j$, we assume that we
have strict inequalities in \eqref{ec.implicit}, then  $(Z^j)$
verifies

\begin{equation}
\label{supersol}
\begin{array}{lcl}
 M \partial Z^{j+1} & > & -A Z^{j+1} +
M((\overline{U}^j)^p - (\underline{U}^j)^p), \\
Z^0 & > & 0. \end{array}
\end{equation}

If the statement of the Lemma is false, then there exists a first
time $t^{j+1}$ and a node $x_i$ such that $z_i^{j+1} \le 0$. There
we have

$$
z_i^{j+1} > z_i^j - \tau_j \frac{a_{ii}}{m_i} z_i^{j+1}  -
\sum_{k\ne i} \frac{a_{ik}}{m_i}z_k^{j+1} + (\overline{u}_i^j)^p -
(\underline{u}_i^j)^p \ge 0,$$ a contradiction.

\subsection{When does the numerical solution blow up}

%
%
%

\begin{lema}
\label{siesgrandeexplota.i} If  $(U^j)_{j\ge 0}$ is large enough,
then blows up in finite time. Moreover
$$\|U^j\|_\infty \sim w^j \sim j$$
\end{lema}

\dem
\begin{eqnarray}
\nonumber w^{j+1} & = & w^j - \tau_j \sum_{k=1}^N\sum_{i=1}^N
a_{ik}u_k^{j+1} + \tau_j
\sum_{i=1}^{N}m_i (u_i^j)^p\\
\nonumber & \le & w^j  + \tau_j \sum_{i=1}^{N} m_i (u_i^j)^p\\
\nonumber & \le  &  w^j   + C \tau_j (w^j)^p\\
\nonumber & = & w^j + C\lambda.\\
\nonumber
\end{eqnarray}

Hence $w^j \le Cj$. To prove the inverse inequality we observe
that

\begin{eqnarray}
\nonumber w^{j+1} & = & w^j - \tau_j \sum_{k=1}^N\sum_{i=1}^N
a_{ik}u_k^{j+1} + \tau_j
\sum_{i=1}^{N}m_i (u_i^j)^p\\
\nonumber & \ge & w^j - \tau_j C_1 w^{j+1} + \tau_j C_2 (w^j)^p\\
\end{eqnarray}
that is

\begin{equation}
\label{pepe}
 (1+C_1\tau_j)w^{j+1} \ge w^j + C_2 \tau_j (w^j)^p .
\end{equation}
Now we look for a subsolution of \eqref{pepe} of the form
$z^j=\Gamma j$. This sequence verifies
$$(1+C_1\tau_j)z^{j+1}= z^j + \Gamma C_1 \tau_j j +
\Gamma(1+C_1\tau_j) \le z^j + C_2 \tau_j (z^j)^p$$ if $\Gamma$ is
small enough. As the discrete maximum principle holds for this
equation we obtain

$$w^j \ge z^j = \Gamma j.$$
This completes the proof. \qed

\medskip

\begin{teo}
\label{funcional.i} Positive solutions of \eqref{td} blow up if
there exists $j_0$ such that $\Phi_{h}(U^{j_0})<0$.
\end{teo}

\dem Again we first observe that $\Phi_{h} (U^j)$ decreases with
$j$, to do that we take inner product of \eqref{ec.implicit} with
$U^{j+1}-U^j$ to obtain

\begin{eqnarray}
\nonumber  0 & = & \langle \frac{1}{\tau_j} M(U^{j+1} -  U^j)  +
AU^{j+1}
- M(U^j)^p,  U^{j+1} -  U^j\rangle\\
\nonumber &=&\tau_j \langle M\partial U^{j+1}, \partial
U^{j+1} \rangle + \Phi_h(U^{j+1}) - \Phi_h(U^j) + \frac12 \langle AU^{j+1}, U^{j+1} \rangle\\
\nonumber & & -
  \langle AU^j,
U^{j+1}\rangle + \frac12 \langle AU^j, U^j\rangle + \frac{p}{2
}\langle M(\xi^j)^{p-1}, (U^{j+1} -  U^j)^2 \rangle.
\end{eqnarray}
Hence we obtain,
\begin{eqnarray}
\nonumber \Phi_h(U^{j+1}) - \Phi_h(U^j) & = & -\tau_j \langle
M\partial U^{j+1},
\partial U^{j+1} \rangle - \frac{\tau_j^2}{2}\langle A\partial
U^{j+1}, \partial U^{j+1} \rangle \\
\nonumber & &- \frac{p}{2
}\langle M(\xi^j)^{p-1}, (U^{j+1} -  U^j)^2 \rangle \\
\nonumber & \le & 0.
\end{eqnarray}
This implies that $\Phi_h$ is a Lyapunov functional for
\eqref{ec.implicit}. We observe that the steady states of
\eqref{ec.implicit} are the same of \eqref{ec.explicit}, so they
have positive energy. Now, assume $(U^j)$ is a bounded solution of
\eqref{ec.implicit}, then there exists a convergent subsequence of
$(U^j)$ that we still denote $(U^j)$. Its limit $W$ is a steady
state with positive energy.

As $\Phi_{h} (U^j)$ decreases and there exists $j_0$ such that
$\Phi_{h} (U^{j_0}) < 0$ then $\Phi_{h,\lambda} (W) < 0$, a
contradiction. We conclude that $(U^j)$ is unbounded and by Lemma
\ref{siesgrandeexplota.i} has finite time blow-up. \qed

\begin{cor}
Assume the convergence hypotheses {\rm (H1), (H2)}. Let $u_0$ an
initial data for \eqref{1.1} such that $u$ blows up in finite time
$T$. Then $u_{h, \lambda}$ also blows up in finite time $T_{h,
\lambda}$ for every $h$, $\lambda=\lambda(h)$ small enough.
Moreover
$$
\lim_{h\to 0} \lim_{\lambda \to 0} T_{h, \lambda} = T.
$$
\end{cor}

\dem If $u$ blows up in finite time $T$ then (see
\cite{CPE},\cite{GK2})

$$ \Phi (u) (t) \equiv \int_\Omega \frac{|\nabla u (s,t) |^2  }{2} \
ds-
 \int_\Omega \frac{ (u(s,t))^{p+1} }{p+1} \ ds \to -\infty \qquad (t
 \nearrow T).$$
Hence there exists a time $t_0<T$ with $\Phi(u)(t_0)<0$. Now we
use the convergence hypothesis (H1) and the convergence of
$u_{h,\lambda}$ to $u_h$ in $[0,t_0]$ to see that
$$
\lim_{h\to 0}\lim_{\lambda \to 0}
\Phi_{h,\lambda}(u_{h,\lambda})(t_0) = \Phi(u)(t_0).
$$
So for $h,\lambda(h)$ small enough we get
$\Phi_{h,\lambda}(u_{h,\lambda})(t_0)<0$ and so $(U^j)$ blows up.
The convergence of the blow-up times is obtained like in the
explicit scheme.

\qed

Next we turn our attention to the  blow-up rate of the discrete
solutions.

\subsection{Blow-up rate}

\begin{teo}
\label{tasa.i}
 Let $u_{h, \lambda}$ a solution with blow-up at time
$T_{h,\lambda}$, then

$$\displaystyle \max_{1 \le i \le N} u_i^j \sim (T_{h,\lambda} -
t^j)^{-1/(p-1)}.$$ Moreover
$$ \displaystyle \lim _{j \to \infty} \max_{1 \le i \le N} u_i^j
(T_{h,\lambda} - t^j)^{1/(p-1)} = C_p = \left (
\frac{1}{p-1}\right )^{1/(p-1)}.$$

\end{teo}

\dem The first part of the proof is the same as the one for the
explicit scheme so we assume we have proved

$$\|U^j\|_\infty \sim (T_{h,\lambda} - t^j)^{-\frac{1}{p-1}},$$
and we are going to recover the constant $C_p$. We change
variables as in the explicit scheme, let $(Y^j)$ be defined by

\begin{equation}
y_i^j=u_i^j(T_{h,\lambda} - t^j)^{1/(p-1)} \qquad 1 \le i \le N.
\end{equation}

In the sequel of the proof we will use $\Delta y_i^{j+1}$ to
denote

$$\frac{y_i^{j+1} - y_i^j}{\tau_j/(T_{h, \lambda}-t^j)},$$
this can be thought as $\tau_j/(T_{h,\lambda})$ to be the time
step in the new variables. With this notation the new variables
verify

\begin{eqnarray}
\label{autosimilares}\nonumber m_i\Delta y_i^{j+1} &=&-(T_{h,
\lambda}-t^j) \displaystyle \sum_{i=1}^N a_{ki} y_{i}^{j+1} +
m_i\frac{(T_{h, \lambda}-t^{j+1})^{1/(p-1)}}{(T_{h,
\lambda}-t^j)^{1/(p-1)}}(y_i^j)^p
\\
\nonumber & - & m_iu_i^j((T_{h, \lambda}-t^{j})^{1/(p-1)}-(T_{h,
\lambda}-t^{j+1})^{1/(p-1)}),\\
\\
\nonumber y^0_i &=&  T_{h, \lambda}^{1/(p-1)}u_0(x_i).\\
\end{eqnarray}

If we assume that there exists $\ve >0$ and a subsequence such
that $y_i^j > C_p + \ve$ for some $i=i(j)$. Then for those
$y_i^j$,  as they are bounded, if $j$ is large enough we have

\begin{eqnarray}
\label{Cp.i}
 \nonumber \Delta y_i^{j+1} &\ge & -\delta +
m_i\frac{(T_{h, \lambda}-t^{j+1})^{1/(p-1)}}{(T_{h,
\lambda}-t^j)^{1/(p-1)}}\left (
(y_i^j)^p - \frac{1}{p-1}y_i^j \right ) \\
\nonumber &  & + \frac{1}{p-1}y_i^j \left[ 1 - \frac{(T_{h,
\lambda}-t^{j+1})^{\frac{1}{p-1}-1}}{(T_{h,
\lambda}-t^j)^{\frac{1}{p-1}-1}} \right ]\\
\nonumber & \ge & \delta.\\
\end{eqnarray}

This means that actually $y_i^j > C_p + \ve$ for every $j$ large
enough and consequently \eqref{Cp.i} is verified for all those
$j$. So $y_i^j$ is unbounded, a contradiction.

The case where there exists an infinite number of $y_i^j$ such
that $y_i^j < C_p - \ve$ is very similar. So we can conclude that
as $j \to \infty$ either $y_i^j \to 0$ or $y_i^j \to C_p$. Now we
use the blow-up rate to obtain
$$\lim_{j \to \infty}\max_{1 \le i \le N} y_i^j =
 \lim_{j \to \infty}\max_{1 \le i \le N} (T_{h, \lambda} -
 t^j)^{1/(p-1)}u_i^j= C_p,
$$
as we wanted to prove. \qed

\subsection{Blow-up set} The propagation property for the blow-up
nodes holds for the implicit scheme and its proof is very similar.
We do not include it.

\section[]{Appendix}
\setcounter{equation}{0}
 In this appendix we prove that if the
general method considered for the space discretization is
consistent (see below) then the totally discrete method converges
in the $L^\infty$ norm. We perform the proofs for the explicit
scheme, they can be extended to the implicit one.

\begin{defi}
\label{consistency}
Let $w$ be a regular solution of
$$
\begin{array}{ll}
w_t = \Delta w + f(x,t) & \qquad \mbox{in } \Omega\times(0,T),\\
w = 0 &  \qquad \mbox{on } \partial\Omega\times(0,T).
\end{array}
$$
We say that the scheme \eqref{semidiscrete} is consistent if for
any $t\in (0,T-\tau)$ it holds
\begin{equation}
\label{assumption} m_i w_t(x_i,t) = -\sum_{k=1}^N a_{ik} w(x_k,t)
+ m_i f(x_i,t) + \rho_{i,h}(t),
\end{equation}
and there exists a function $\rho:\RR_+\to \RR_+$ such that
$$\max_i \frac{|\rho_{i,h}(t)|}{m_i} \le \rho(h), \qquad
\mbox{for every } t\in (0,T-\tau),$$ with $\rho(h)\to 0$ if $h\to
0$. The function $\rho$ is called the modulus of consistency of
the method.
\end{defi}

If we consider for example a finite differences scheme in a cube
$\Omega=(0,1)^d \subset \mathbb{R}^d$. Then the modulus of
consistency can be taken as $\rho(h)=Ch^2$.

\begin{teo}
\label{conver} Let $u$ be a regular solution of {\rm (\ref{1.1})}
($u \in C^{2,1} (\overline{\Omega} \times [0,T-\tau]$) and
$u_{h,\lambda}$ the numerical approximation given by {\rm
(\ref{td})} then there exists a constant $C$ depending on $\|u\|$
in $C^{2,1} (\overline{\Omega} \times [0,T-\tau])$ such that
$$\| u - u_{h,\lambda} \|_{L^\infty ([0,T-\tau], L^2(\Omega))}
\le C (\rho(h) + \lambda). $$
\end{teo}

\dem

We define the error functions
$$
e_i^j = u(x_i,t_j) - u_i^j.
$$
By \eqref{assumption}, these functions verify
$$
m_i \partial  e_i^{j+1} = -\sum_{k=1}^N a_{ik} e_k^j + m_i(
u^p(x_i,t_j)- (u_i^j)^p) + \rho_i(h) + Cm_i\lambda,
$$
where $C$ is a bound for $\|u_{tt}\|_{L^\infty(\Omega \times
[0,T-\tau])}.$  Let
$$t_0 = \max \{t: \ t < T- \tau, \ \max_i \max_{t_j  < t}
|e_i^j|\le 1 \}.$$ We will see by the end of the proof that $t_0 =
T-\tau$ for $h, \lambda$ small enough. In $[0,t_0]$ we have
$$
m_i \partial  e_i^{j+1} = -\sum_{k=1}^N a_{ik} e_k^j +
m_ip(\xi_i^j)^{p-1}e_i^j + \rho_i(h) + Cm_i\lambda,
$$
hence, in $[0, t_0]$, $E^j = (e^j_1,...,e^j_N)$ satisfies
\begin{equation}
\label{eq.error} M \partial E^{j+1} \le -AE^j + KME^j + (\rho(h) +
C\lambda)M1^t.
\end{equation}
Let us now define $W^j=(w^j_1, \dots, w_N(t))$, which will be used
as a supersolution.
$$
w_i^j= e^{(2K+1)t_j}(\|e(0)\|_{\infty} + \rho (h) + C\lambda).
$$
It is easy to check that $W^j$ verifies
$$
M \partial W^{j+1} > -AW^j + KMW^j + (\rho(h) + C\lambda)M1^t,
$$
here $K$ is the Lipchitz constant for $f(u)=u^p$ in $[0,
\|u(\cdot,T-\tau)\|_{L^\infty}]$. Hence $W^j$ is a supersolution
for \eqref{eq.error}, and by Lemma \ref{pmaximo} we get
$$
 e^j_i \le e^{(2K+1)t_j}(\|e^0\|_{L^\infty(\Omega)} + \rho
(h) + C\lambda ), \qquad t_j \in [0, t_0].
$$
Arguing along the same lines with $-E^j$, we obtain
$$
|e^j_i| \le e^{(2K+1)T} (\|E(0)\|_{\infty} + \rho (h) + C\lambda)
\le C( \rho(h) + \lambda), \qquad t_j \in [0, t_0],
$$
by our hypotheses on the convergence of the initial data. Using
this fact, since $\rho(h)$ goes to zero, we get that $|e^j_i| \le
1$ for every $t_j \in [0,T-\tau]$ for every $h, \lambda$ small
enough. Therefore $t_0 = T-\tau$ for $h, \lambda$ small enough.
This proves the convergence of the scheme. In fact we have that
for every  $h <h_0, \lambda < \lambda_0$

$$\max_{j} \max_{1\le i \le N} |u^j_i - u(x_i,t_j)|
\le C(\rho(h) + \lambda).$$ \qed

{\centerline {\bf Acknowledgments}}

I would like to thank Julio D. Rossi for his help, his
encouragement, and his friendship. I really enjoy and learn a lot
working with him.

\end{document}